\documentclass[5p]{elsarticle}

\usepackage{amsmath,amssymb,subfigure}
\usepackage{pstricks,pst-plot}
\usepackage{psfrag}
\usepackage{xcolor}
\newtheorem{theorem}{Theorem}[section]

\journal{Physica D}
\def\bv#1{\mbox{\bfseries\itshape#1}}
\def\Im{\mathrm{Im}\,}
\def\Re{\mathrm{Re}}
\def\imagunit{\sqrt{-1}}
\def\uu{\mbox{\bfseries\itshape u}}

\definecolor{myred}{rgb}{0.8,0.2,0.2}
\definecolor{myblue}{rgb}{0.4,0.5,0.8}
\newgray{mygray}{0.9}

\begin{document}

\begin{frontmatter}


\title{Dynamics of poles with position-dependent strengths\\ and its optical analogues}

\author[authorlabel1]{James Montaldi}
\ead{j.montaldi@manchester.ac.uk}
\author[authorlabel2]{Tadashi Tokieda}
\ead{tokieda@dpmms.cam.ac.uk}

\address[authorlabel1]{School of Mathematics, University of Manchester, Manchester M13 9PL, England}
\address[authorlabel2]{Trinity Hall, Cambridge CB2 1TJ, England}

\begin{abstract}

The dynamics of point vortices is generalized in two ways: first by making the strengths complex, which allows for sources and sinks in superposition with the usual vortices, second by making them
functions of position.  These generalizations lead to a rich dynamical system, which is nonlinear 
and yet has conservation laws coming from a Hamiltonian-like formalism.  We then discover 
that in this system the motion of a pair mimics the behavior of rays in geometric optics.  We describe several exact solutions with optical analogues, notably Snell's law and the law of reflection off a mirror, and perform numerical experiments illustrating some striking behavior.
\end{abstract}


\end{frontmatter}

\section*{Introduction}

   The dynamics of point vortices with fixed strengths in a 2-dimensional ideal fluid has a classical pedigree, e.g.\ Lamb \cite{Lamb}, Art.\ 154--160.   The subject continues to be actively pursued: a modern survey \cite{Aref} on its equilibrium aspect alone lists more than 100 papers.  We generalize vortex dynamics in two ways, firstly allowing, besides vortices, sources/sinks as well as their superpositions (`poles'), and secondly allowing the strengths of these poles to vary as functions of position in the plane.

The first generalization goes back to a 1928 paper by Friedman and Polubarinova-Kochina (the former is the same Friedman as in the eponymous cosmological model).  The rather more recent paper by Borisov \&\ Mamaev \cite{BM} contains references as well as a good theoretical analysis; see also a paper by Lacomba \cite{Lacomba}.  Here we present a couple of new exact solutions and alternative derivations of some old ones.  

The second generalization seems less explored, and leads to rich dynamics, which we illustrate with a variety of exact solutions with {\it analogues in geometric optics}, the position-dependent strength of a pole replacing the medium-dependent index of refraction. As typical examples we detail the analogues of Snell's law and the law of reflection off a mirror, in generalized forms.  Optical analogy is not so obvious: though it was suggested by Kimura \cite{Kimura} that in the dipole limit a classical vortex pair should travel along a geodesic, the principle governing light rays in geometric optics is one of least time, not of least length.  Nevertheless, the dynamics of poles with position-dependent strengths turns out to be quite versatile in mimicking the ray representation of phenomena of wave propagation.  Take for instance the work by Berry \cite{Berry} on focusing and defocusing of surface waves by underwater landscape.  It will become clear that such effects are realizable by our dynamics, too.  There are also similarities with results of Longuet-Higgins on trapping waves around islands \cite{LH}.

In an interesting paper, Hinds, Johnson and McDonald \cite{HJM} consider the dynamics of a pair of vortex patches as they cross a step change in the depth of the fluid, and also find the pair is refracted provided they are sufficently well separated compared to their size and the angle of incidence, otherwise they find vortex shedding.  While we do not claim our introduction of a `seabed' function genuinely models a variation in depth, we do consider only point vortices (and poles) so vortex shedding would not arise.

In the language of dynamical systems, this second generalization through the introduction of a `seabed' function $S$ in section 4 puts us in the realm of {\it hybrid systems}, where different equations of motion govern different regions of the phase space. 

   It may not be amiss to point out that the dynamical system (2), which is the chief object of our study, is quite nonlinear---in 
a sense more so than say the Euler or Navier-Stokes equations.  In the latter, indeed, the nonlinearity is separated out as an additive term $(\bv{v}\cdot \nabla)\bv{v} = \nabla (\frac{1}{2}\bv{v}^2) - \bv{v}\times (\nabla \times \bv{v})$, so that we can resort to linearization by dropping this nonlinear term 
or by substituting for it a term $({\rm background~flow}\cdot \nabla)\bv{v}$ (Stokes and Oseen approximations, respectively).  In contrast, our system is so nonlinear that it is not even clear whether or
not there exists a `linearization' of the system that makes sense.  

   In the light of this inseparable nonlinearity, it is noteworthy that our system proves to be analysable in elementary
terms, and many instances of complete integrability explicitly spelt out.  This is thanks, 
ultimately, to the fact that the classical vortex dynamics is {\it Hamiltonian\/} and our generalization is
{\it something like Hamiltonian}.

\section{Equation of motion}

We begin with a discussion of 2-dimensional ideal fluid flow in terms of complex potential, but shall specialize soon.   Consider $N$ interacting points $z_1, \ldots, z_N\in {\mathbb C}$ called {\bf poles}, each pole $z_i$ carrying with it a family of complex-valued functions $\bigl\{ \mu_n^i(z) \bigr\}_{n\in {\mathbb Z}}$ called {\bf strengths}, only finitely many of which are nonzero.  The poles move according to
\begin{equation}\label{eqn of motion}
\frac{d}{dt} z_i(t) = \sum_{j\, :\, j\neq i}\, \sum_{n\in {\mathbb Z}}\,
\overline{ \mu_n^j\bigl(z_j(t)\bigr)\, \bigl(z_i(t) - z_j(t)\bigr)^{\! n} }
\end{equation}
for $i = 1, \ldots, N$
(${}^{\overline{ {\quad} }}$ denotes complex conjugation).  The dynamics of (\ref{eqn of motion}), being 1st-order in time $t$, has {\it no inertia}, in the sense that the poles' instantaneous positions determine their velocities: the phase space is a product of $N$ copies of $\mathbb C$ (minus diagonals if we wish a priori to exclude collisions), not a (co)tangent bundle.  We can set up a dynamical system like this on any domain of any Riemann surface.
In simple domains that arise in practice, solutions can be found by the method of images. 

A term of the form $\overline{ \mu (z - z_j)^n }$ on the right side of (\ref{eqn of motion}) represents a flow velocity induced by $z_j$ at $z$.  The pictures for $n = -1$ have rotational symmetry: source or sink of flux $2\pi \mu$ for $\mu$ real, vortex of circulation $2\pi \imagunit \, \mu$ for $\mu$ pure imaginary, in general a superposition of these, i.e.\ a spiral node.  Other $n$ exhibit other symmetries: multipolar flows for $n < -1$, uniform flow for $n = 0$, and corner flows for $n > 0$.


\section{Homogeneous systems, conservation laws}

Now suppose (\ref{eqn of motion}) is homogeneous so that \ $\mu_n^i = 0$ except for a certain exponent $n = n_0$, and moreover all $\mu^i := \mu_{n_0}^i$ are fixed.  Then (\ref{eqn of motion}) may be recast in the `canonical' form
$$
\frac{d}{dt} z_i = \frac{2}{\overline{\mu^i}} \frac{\partial}{\partial \, \overline{z_i}} H\, ,
$$
with
$$
H(z_1, \ldots, z_N) = \Re \!\sum_{i, j\, :\,  i < j}\,  \mu^i\mu^j\, G(z_i - z_j)
$$
and
\begin{equation*}
G(z_i - z_j) = 
\begin{cases} 
\frac{1}{n_0 + 1} (z_i - z_j)^{n_0 + 1} & \text{when $n_0 \neq -1$}\, , \\[6pt]
\log (z_i - z_j) & \text{when $n_0 = -1$}
\end{cases}
\end{equation*}
($G$ as in `Green').  From
\begin{eqnarray*}
\frac{d}{dt}H\bigl(z_1(t), \ldots, z_N(t)\bigr) \hskip -15mm && \\
&=& 
\sum_i \Bigl( \frac{\partial}{\partial z_i} H \cdot \frac{d}{dt} z_i +
 \frac{\partial}{\partial \, \overline{z_i}} H \cdot \frac{d}{dt} \overline{z_i} \Bigr) \\
&=& \sum_i \Re(\mu)\, \left| \frac{d}{dt}z_i(t) \right|^2
\end{eqnarray*}
we see

\smallskip

\begin{theorem}\label{thm:conservation of H}
If all $\mu^i$ are pure imaginary and fixed, then $H$ is conserved.
\end{theorem}

\smallskip

   Next let the homogeneity degree $n_0$ be odd, with $\mu^i$ still fixed.  Pairwise
cancellation in (\ref{eqn of motion}) yields $\sum_i \overline{\mu^i}\, dz_i/dt = 0$, whence

\smallskip

\begin{theorem}\label{thm:centre of strength}
If the degree of homogenity is odd and  $\mu := \sum_i \mu^i \neq 0$, then the `center
of strength'  
$$
c = \frac{\sum_i \overline{\mu^i}\, z_i}{\overline{\mu}}
$$ 
is conserved.
If $\mu = 0$, then for every partition of the index set
$I\sqcup I' = \{ 1, \ldots, N \}$ such that $\mu_I := \sum_{i\in I} \mu^i \neq 0$ and
$\mu_{I'} := \sum_{i'\in I'} \mu^{i'} \neq 0$, the difference between the `subcenters'
$$
\frac{\sum_{i\in I} \overline{\mu^i}\, z_i}{\overline{\mu_I}} -
\frac{\sum_{i'\in I'} \overline{\mu^{i'}}\, z_{i'}}{\overline{\mu_{I'}}}
$$ 
is conserved.
\end{theorem}

\smallskip

\noindent The `partition' part of this Theorem is elementary but does not appear to have been used prior to the paper by Montaldi, Souli\`ere, Tokieda \cite{MST}.
  
\medskip

If instead the homogeneity degree $n_0$ is even and there are just 2 poles $z,z'$ with strengths $\mu,\mu'$, then 
$\overline{\mu}z-\overline{\mu'}z'$ is a conserved quantity.  However, this does not appear to extend to more than 2 poles.

\medskip

Finally, how can we extend the affine symmetry of the phase space to that of the phase space-time so as to preserve the invariance of (\ref{eqn of motion})?  The requirement that time $t$ be real gives the answer.

\smallskip

\begin{theorem}\label{thm:invariance}
The system (\ref{eqn of motion}) is invariant under the action of ${\mathbb C}^* \ltimes {\mathbb C}$ if and only if it is homogeneous of degree $n_0 = -1\,$: here $(a, b)\in {\mathbb C}^*\ltimes {\mathbb C}$ acts by sending $(t, z)$ to $(|a|^2t, az + b)$.
\end{theorem}

\section{Exact solutions}

   We shall, in the remainder of the paper, focus on the theory where (\ref{eqn of motion}) involves only the exponent $n_0 = -1$.  In this section we suppose all the strengths are fixed: $\mu^i_{-1}(z)=\mu^i$.  Thus the equations of motion become
\begin{equation}\label{eq:dz_i/dt fixed mu}
\frac{d}{dt} z_i(t) = \sum_{j\, :\, j\neq i} \frac{ \overline{\mu^j} }{ \overline{z_i(t)} - \overline{z_j(t)} }
\qquad (i = 1, \ldots, N)\, .
\end{equation}
We consider position-dependent strengths in sections 4 and 5.  When  all $\mu^i$ are pure imaginary, we are back to classical point vortices and recover the logarithmic $H$ as their Hamiltonian.

\subsection{Self-similar solutions}

If a collection of poles happens to move in a self-similar manner, then Theorem \ref{thm:invariance} reduces the (complex) degree of freedom from $N$ to $1$, down to a single equation
\begin{equation}\label{eq:dZ/dt}
\frac{d}{dt} Z = \frac{M}{\, \overline{Z}\, }\, ,
\end{equation}
or in polar coordinates
$$
\frac{d}{dt} \frac{1}{2}|Z|^2 = {\rm Re}\, M\, , \qquad |Z|^2\frac{d}{dt} \arg Z = {\rm Im}\, M\, .
$$
The solution is
\begin{equation}\label{eq:Z(t)}
Z(t) = T \exp\Bigl( \imagunit \, \frac{{\rm Im}\, M}{{\rm Re}\, M}\, \log T \Bigr) Z(0)\, ,
\end{equation}
where 
$$T = \sqrt{1 + \frac{{\rm Re}\, M}{|Z(0)|^2/2}\, t\,}
\qquad {\rm if~~Re}\, M \neq 0\, ,$$
and
\begin{equation}\label{eq:Z(t)M=0}
Z(t) = \exp\Bigl( \imagunit \, \frac{{\rm Im}\, M}{|Z(0)|^2}\, t \Bigr) Z(0)
\qquad 
{\rm if~~Re}\, M = 0\, .
\end{equation}
We remark however that there is really no need to treat the case (\ref{eq:Z(t)M=0}) apart from (\ref{eq:Z(t)}), since (\ref{eq:Z(t)}) converges to (\ref{eq:Z(t)M=0}) in the limit ${\rm Re}\, M \to 0$.

   If ${\rm Re}\, M \neq 0$, then the poles spiral in and collapse to their center of strength after time
$- |Z(0)|^2/2\,{\rm Re}\, M$  (in the future if ${\rm Re}\, M < 0$, in the past if $ > 0$).  If ${\rm Re}\, M = 0$, then the configuration of the poles merely spins while remaining congruent to itself.

\subsection{Pair} \label{sec:pair}

   A pole pair $z$, $z'$ of fixed strengths $\mu$, $\mu'$ moves self-similarly around its center of strength $c$
.  With a little manipulation (\ref{eq:dz_i/dt fixed mu}) takes the form of (\ref{eq:dZ/dt}):
$$
\frac{d}{dt} (z - c) = \frac{|\mu'|^2}{\mu + \mu'} \frac{1}{ \overline{z} - \overline{c} }
$$
and the same equation for $z' - c$ with prime and unprime exchanged.  
If $\mu + \mu' = 0\,$, then Theorem 2.2 implies that
the pair moves along parallel trajectories, around a `center at infinity' ($c \to \infty$).   In particular if $\mu$, $\mu'$ are real and $\mu + \mu' = 0$, then the source chases the sink and the sink runs away from the source, while their mutual distance remains constant.

\medskip

More generally, a short calulation shows that the distance between a pair of poles remains constant provided their strengths have the opposite real part: $\Re(\mu+\mu')=0$. In other words, one should be a sink, the other a source of the same strength, and each combined with an arbitrary rotational (vortex) part. Indeed, we find from (\ref{eq:dz_i/dt fixed mu}) that
$$\frac{d}{dt}|z-z'|^2 = 2\Re(\mu+\mu').$$

\subsection{Regular polygons}

   Place $N$ poles of equal fixed strengths $\mu$ at the vertices of a regular $N$-gon,
plus 1 pole of strength $\mu'$ at the center $c\,$.  The poles spiral self-similarly out or in,
clockwise or anticlockwise, depending on $N, \mu, \mu'$.   In terms of any one of the vertices $z\,$,
(\ref{eq:dz_i/dt fixed mu}) takes the form of (\ref{eq:dZ/dt}):
$$
\frac{d}{dt} (z - c) = \left( \frac{N-1}{2}\, \overline{\mu} + \overline{\mu'} \right) \frac{1}{ \overline{z} - \overline{c} }\, .
$$
If $\mu' = - (N - 1)\mu/2$, then all the poles are immobilized.

\section{Position-dependent strengths, optical analogy}

We continue to assume that the equation of motion (\ref{eqn of motion}) involves only the exponent $n_0 = -1$ but now we allow the strengths to be position-dependent, in such a way that they all share a common dependence on the position: $\mu^i(z) = \mu^i\cdot S(z)$ for some complex constants $\mu^i$ $(i = 1,\ldots, N)$ and some real-valued function $S$ (`seabed' or `step').  Thus the equations of motion are
 \begin{equation}\label{eq:dz_i/dt}
 \frac{d}{dt} z_i(t) = \sum_{j\, :\, j\neq i} \frac{ \overline{\mu^j} S(z_j(t))}{ \overline{z_i(t)} - \overline{z_j(t)} }
 \qquad (i = 1, \ldots, N)\, .
 \end{equation}

We begin by deriving a Noether-type conservation law under the extra assumption that all $\mu^i$ are imaginary (vortices); on the other hand, the argument that follows works for the general exponent $n_0$ . Now under this extra assumption, the differential 2-form
$$\Omega = \frac12\sum_i \mu^iS(z_i)\, dz_i\wedge d\overline{z_i}.$$
is a \emph{symplectic form} on the phase space: it is closed $d\Omega = 0$ and non-degenerate wherever all $S(z_i)$ are nonzero, though the non-degeneracy is not needed for the argument. Denote the vector field (\ref{eq:dz_i/dt}) by $V$, and let  
\begin{equation}\label{eq:hamiltonian}
H(z_1, \ldots, z_N) = \Re \!\sum_{i, j\, :\,  i < j}\,  \mu^i\mu^j\,S(z_i)S(z_j) G(z_i - z_j).
\end{equation}
Unless the strengths are constants, Hamilton's equation does not hold: $dH(\cdot)$ is not equal to $\iota_V\Omega(\cdot) = \Omega(V, \cdot )$,
because of the terms involving the derivative of $S(z)$. However, suppose $S$ admits the symmetry of some Euclidean motion represented by a vector field 
$\uu$ on $\mathbb{C}$---meaning that \uu\ is everywhere tangent to the level curves of $S$ and preserves mutual distances between points. Then the \uu-component of Hamilton's equation does hold, and we have:
$$\iota_{\uu}\Omega(V)=\Omega(\uu,V)=-\Omega(V,\uu) = -dH(\uu)=0.$$
Abusing notation, we denote by \uu\ both the given vector field on $\mathbb C$ and its induced vector field on the phase space $\mathbb{C}^N$, 
$\uu(z_1,\dots,z_N) = (\uu(z_1),\dots,\uu(z_N)).$
The \emph{momentum conjugate} to \uu\ is a function $\psi$ on the phase space such that $d\psi = \iota_{\uu}\Omega$ ($\psi$ is just the stream function for \uu). The above equation says that $d\psi(V ) = 0$. Thus,

\begin{theorem}\label{thm:Noether}
If $S$ admits the symmetry of a Euclidean motion \uu, then the momentum $\psi$ conjugate to \uu\ is a conserved quantity of the dynamics.
\end{theorem}

There are two types of Euclidean motion, rotation and translation. We now write $\psi$ explicitly for both types. For rotation, $\uu = (-y, x)$ in Cartesian coordinates, we have 
$$\iota_{\uu}\Omega = -\sum_i\imagunit\,\mu^i S(z_i)(x_i\,dx_i + y_i\,dy_i)$$
(since $dz\wedge d\overline{z}=-2\imagunit\, dx\wedge dy$). If $S$ admits a rotational symmetry about the origin, we can write $S = S(|z|^2 )$, and then we find $d\psi = \iota_{\uu}\Omega$ is satisfied by

\begin{equation}\label{eq:rotational momentum}
\psi(z_1,\dots,z_N) = -\sum_i \imagunit\,\mu_i\sigma(|z_i|^2),
\end{equation}
where $\sigma'(u)=S(u)$.  For translation, say $\uu = (1, 0)$, we find
\begin{equation}\label{eq:translational momentum}
\psi(z_1,\dots,z_N) = -\sum_i \imagunit\,\mu_i\sigma(y_i)
\end{equation}
 where $\sigma'(y)=S(y)$. If $S$ is constant, then these become the usual conserved quantities in the point vortex problem: respectively the angular impulse and (one component of) the linear impulse, or center of vorticity after rescaling.

\medskip

In the case of a pair of vortices, i.e.\ $\mu$, $\mu'$ are pure imaginary, their mutual distance remains constant, even with an arbitray `seabed' function $S(z)$.  This is easy to see, for the velocities of the vortices are then perpendicular to the segment connecting the vortices; whether the velocities have the same magnitude or not is irrelevant.

In most of our examples $S(z)$ will be piecewise constant.  Then Theorems 2.1, 2.2, 2.3 hold piecewise until one of the poles crosses a discontinuity of $S$, at which instant $H$ and the centers of strength jump to new values. Moreover, in each of the examples $S$ has a symmetry 
for which our Noether-type theorem gives rise to a conservation law.

\subsection{Analogue of Snell's law}\label{sec:refraction}

   Let $S(z) = s_1$ in the lower half-plane ${\rm Im}\, z < 0$ and $S(z) = s_2$ in the upper half-plane ${\rm Im}\, z \geqslant 0$, where
$s_1, s_2\in {\mathbb R}_+$.  A pair of poles $z$, $z'$ with strengths $- \mu S(z)$,
$+\mu S(z')$ ($\mu\in\mathbb{C}$), moving along parallel trajectories, arrives from the lower half-plane (Figure~\ref{fig:Snell} shows the case of vortices: $\Re(\mu)=0$).

 \psset{arrowsize=5pt}
 \begin{figure}[htb]
 \psset{unit=1.3}
   \begin{center}
   \begin{pspicture}(-2,-2)(3,2)
   \psline(-2,0)(3,0) 
   \psframe[fillstyle=solid,fillcolor=mygray,linestyle=none](-3,-2.3)(4,0)
 {\psset{linecolor=myred,linewidth=1.5pt}
   \rput(2.5,-0.8){$z'(t)$}
   \psline{->}(2.5125,-1.35)(2.25,-1)
   \psline(2.25,-1)(1.5,0)
   \psarc(-1.5,-2.25){3.75}{36}{65}
   \psline{->}(0.084,1.148)(-1.516,1.975)
}
 {\psset{linecolor=myblue,linewidth=1.5pt}
   \rput(0.9,-1.8){$z(t)$}
   \psline{->}(1.5125,-2.1)(1.25,-1.75)
   \psline(1.25,-1.75)(0.5,-0.75)
   \psarc(-1.5,-2.25){2.5}{36}{65}
   \psline{->}(-0.4,0)(-2,0.7272)
}
   \psline[linestyle=dashed](1.5,0)(0.5,-0.75)
   \psline[linestyle=dashed](-0.4,0)(0.084,1.148)
   \psline[linestyle=dotted](1.5,0)(1.5,-0.75)
   \psline[linestyle=dotted](-0.4,0)(-0.4,0.75)
   \psarc(1.5,0){0.3}{270}{306}
   \rput(1.7,-0.5){$\theta_1$}
   \psarc(-0.4,0){0.3}{90}{155}
   \rput(-0.7,0.4){$\theta_2$}
   \rput(-2.5,0.5){$S=s_2$}
   \rput(-2.5,-0.5){$S=s_1$}
   \end{pspicture}
   \end{center}
   \caption{Snell's law for vortices}\label{fig:Snell}
 \end{figure}
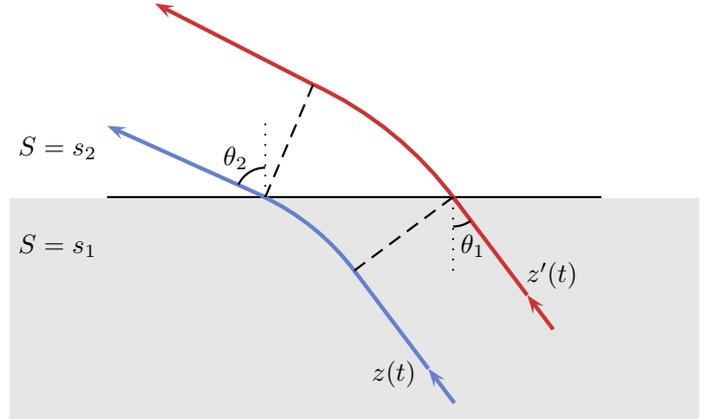
When $z'$ crosses the real axis, the pair starts swerving
self-similarly around its center of strength $c = (s_2 z' - s_1
z)/(s_2 - s_1)$.  When, after time $t$ say, $z$ crosses
the real axis, the pair resumes moving along parallel trajectories inside the upper
half-plane. The angles of incidence $\theta_1$, $\theta_2$ of these trajectories pre- and post-crossing satisfy, if $\mathrm{Im}\, \mu \neq0$,
\begin{equation}\label{eq:generalized Snell}
s_1 \sin\theta_1 \exp\Bigl(-\frac{{\rm Re}\, \mu}{{\rm Im}\, \mu}\, \theta_1 \Bigr) =
s_2 \sin\theta_2 \exp\Bigl(-\frac{{\rm Re}\, \mu}{{\rm Im}\, \mu}\, \theta_2 \Bigr)\,,
\end{equation}
and
\begin{equation}\label{eq:theta1=theta2}
\theta_1 = \theta_2 
\qquad 
{\rm if~~Im}\, \mu = 0\, .
\end{equation}
If $\mu$ is pure imaginary, the case depicted in Figure~\ref{fig:Snell}, then (\ref{eq:generalized Snell}) reduces to
\begin{equation}\label{eq:Snell}
s_1\sin\theta_1 = s_2\sin\theta_2
\end{equation}
which is the analogue of Snell's law in optics; $s$ is
then the analogue of the index of refraction.  At one extreme $\theta_1 = \theta_2 = 0$
there is no refraction.  At the other extreme $\theta_2 = \pi/2$ the pair `skids' along the 
real axis, hence reflection into the lower half-plane occurs for 
$$
\theta_1 > \arcsin \frac{s_2}{s_1} \, .
$$  

\smallskip

   One proof of (\ref{eq:generalized Snell}) and (\ref{eq:theta1=theta2}) goes as follows.  During the crossing, i.e.\ while $z'$ is already in the 
upper half-plane but $z$ is still in the lower half-plane, $z' - z$ evolves according to
$$
\frac{d}{dt} (z' - z) = \frac{ \overline{\mu} (s_2 - s_1)}{ \overline{z'} - \overline{z} }
$$
which is in the form (\ref{eq:dZ/dt}) with $Z = z' - z$ and $M = \overline{\mu}(s_2 - s_1)$.
As remarked after (\ref{eq:Z(t)M=0}), the special case of ${\rm Re}\, \mu = 0$ may be treated as the
limit ${\rm Re}\, \mu \to 0$, so we shall proceed with the proof assuming ${\rm Re}\, \mu \neq 0$.
With this in mind, we have ${\rm Im}\, M/{\rm Re}\, M = - {\rm Im}\, \mu/{\rm Re}\, \mu$, so by (\ref{eq:Z(t)})
\begin{equation}\label{eq:z'-z}
z'(t) - z(t) = T \exp\Bigl( - \imagunit \, \frac{{\rm Im}\, \mu}{{\rm Re}\, \mu}\, \log T \Bigr) (z'(0) - z(0))\, ,
\end{equation}
where
$$
T = \sqrt{1 + \frac{(s_2 - s_1){\rm Re}\, \mu}{|Z(0)|^2/2}\, t\,}\,.
$$
Since
\begin{equation}\label{eq:z-c}
z - c = \frac{s_2}{s_1 - s_2}(z' - z)
\end{equation}
(of course: the point is that during the crossing the solution is self-similar), (\ref{eq:z'-z}) also gives
\begin{equation}\label{eq:z(t)-c}
z(t) - c = T \exp\Bigl( - \imagunit \, \frac{{\rm Im}\, \mu}{{\rm Re}\, \mu}\, \log T \Bigr) (z(0) - c)\, .
\end{equation}
The equations (\ref{eq:z'-z}) and (\ref{eq:z(t)-c}) convey all the information we need to compare the positions of the 
pair at time $0$ (or more precisely $0+$) when the crossing starts and at time $t$ (or $t-$) 
when it ends.  

   First, taking the argument of (\ref{eq:z'-z}) we find
\begin{equation}\label{eq:theta2}
\theta_2 =  - \frac{{\rm Im}\, \mu}{{\rm Re}\, \mu}\, \log T + \theta_1\, .
\end{equation}
Next, taking the imaginary part of (\ref{eq:z(t)-c}) we find
$$
{\rm left~side} = 0 - {\rm Im}\, c = \frac{s_1}{s_2 - s_1}{\rm Im}\, z(0)
$$
because ${\rm Im}\, z(t) = {\rm Im}\, z'(0) = 0$, and we find
$$
{\rm right~side} = \frac{T\,s_2}{s_1 - s_2}|z'(0) - z(0)| \sin \Bigl( - \frac{{\rm Im}\, \mu}{{\rm Re}\, \mu}\, \log T + \theta_1 \Bigr)
$$
because by virtue of (\ref{eq:z-c}) $z(0) - c$ is expressible as 
$$
\frac{s_2}{s_1 - s_2}|z'(0) - z(0)| \exp(\imagunit \, \theta_1)\, .
$$
Using (\ref{eq:theta2}) and equating the two sides,
$$
T s_2 \sin\theta_2 = s_1 \frac{- {\rm Im}\, z(0)}{|z'(0) - z(0)|} = s_1 \sin \theta_1\, .
$$
But again by (\ref{eq:theta2})
$$
{\rm if~~} {\rm Im}\, \mu = 0\, , \quad {\rm then~~} \theta_2 = \theta_1\, ,
$$
and
$$
{\rm if~~} {\rm Im}\, \mu \neq 0\, , \quad {\rm then~~} T = \exp \frac{{\rm Re}\, \mu}{{\rm Im}\, \mu} (\theta_1 - \theta_2) \, .
$$
The formulae (\ref{eq:generalized Snell}) and (\ref{eq:theta1=theta2}) are proved.

\smallskip

   Another, slightly different, formulation of this generalized Snell's law is to say that during the crossing the quantity
$$s \sin \theta\, \exp \left(- \frac{{\rm Re}\, \mu}{{\rm Im}\, \mu}\, \theta \right)
$$ remains constant if ${\rm Im}\, \mu \neq 0$,
and $\theta$ remains constant if ${\rm Im}\, \mu = 0$.  In this formulation, the law clearly holds in the more general set-up
where a real-valued `seabed' function $S(z)$ depends only on ${\rm Im}\, z$ but otherwise varies arbitrarily: 
just divide the plane into thin strips parallel to the real axis and approximate $S(z)$ by a function constant on 
each strip, like a sloping beach.  

\smallskip

When $\mu$ is pure imaginary, a shorter proof of (\ref{eq:Snell}) is available via Noether's theorem \ref{thm:Noether}.  The invariance under translations along the real direction gives rise to the conserved quantity $\mathrm{Re}\,(\mu S(z')z' - \mu S(z)z)$, and hence of $\Im(S(z')z'-S(z)z)$. From the fact already observed that the  distance between the vortices remains constant, (\ref{eq:Snell}) follows immediately by geometry.

\subsection{Analogue of the law of reflection} \label{sec:reflection}

\begin{figure}[htb]
   \begin{center}
\psset{unit=1.7}
   \begin{pspicture}(-1,-2)(3,1.2)
   \psline(-1.5,0)(3.5,0) 
   \psframe[fillstyle=solid,fillcolor=mygray,linestyle=none](-1.5,1)(3.5,0)
{\psset{linecolor=myred,linewidth=1.5pt}
   \rput(2.65,-1.15){$z'(t)$}
   \psline{->}(3,-2)(2.25,-1)
   \psline(2.25,-1)(1.5,0)
   \psarc(1,-0.375){0.625}{36}{144}  \psline{->}(1.1,0.25)(0.95,0.25)
   \psline{->}(0.5,0)(-0.4366,-1.2488)
   \psline(0.5,0)(-1,-2)
}
{\psset{linecolor=myblue,linewidth=1.5pt}
   \rput(1.6,-1.8){$z(t)$}
   \psline{->}(1.4375,-2)(1.25,-1.75)
   \psline(1.25,-1.75)(0.5,-0.75)
   \psarc(1,-0.375){0.625}{-144}{-36} \psline{->}(0.9,-1)(1.05,-1)
   \psline{->}(1.5,-0.75)(0.5625,-2)
}
   \psline[linestyle=dashed](1.5,0)(0.5,-0.75)
   \psline[linestyle=dotted](1.5,0)(1.5,-0.5)
   \psline[linestyle=dashed](0.5,0)(1.5,-0.75)
   \psline[linestyle=dotted](0.5,0)(0.5,0.5)
   \psarc(1.5,0){0.25}{270}{306}
   \rput(1.65,-0.4){$\theta_1$}
   \psarc(0.5,0){0.25}{90}{234}
   \rput(0.15,0.2){$\theta_2$}
   \rput(-1,0.5){$S=s_2$}
   \rput(-1,-0.5){$S=-s_1$}
   \end{pspicture}
   \end{center}
   \caption{Law of reflection for vortices}\label{fig:reflection}
 \end{figure}
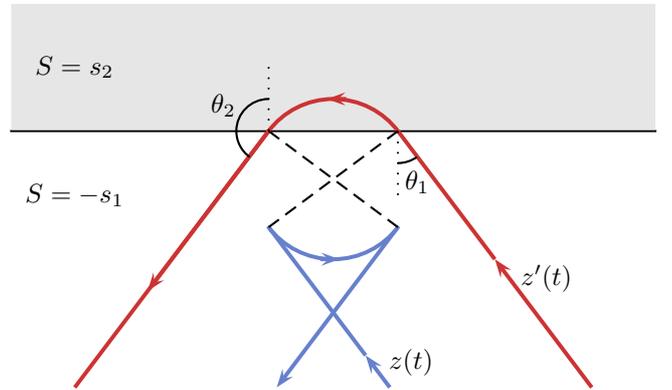

   Let $S(z) = -s_1$ in the lower half-plane ${\rm Im}\, z < 0$ and $S(z) = s_2$ in the upper half-plane ${\rm Im}\, z \geqslant 0$, where $s_1, s_2\in {\mathbb R}_+$.  The difference from the set-up for Snell's law is that $S$ changes sign between the  half-planes.  To fix ideas, let us think of $\mu$ such that ${\rm Re}\, \mu \geqslant 0$, ${\rm Im}\, \mu \leqslant 0$.  
A pair $z$, $z'$ with strengths $- \mu S(z)$, $+\mu S(z')$, moving along parallel trajectories, arrives from the lower half-plane (Figure~\ref{fig:reflection}).

When $z'$ crosses the real axis, the pair starts spiraling out self-similarly around its 
center of strength $c = (s_2 z' + s_1 z)/(s_2 + s_1)$.  When, after time $t$ say, $z'$ again 
crosses the real axis (it does so necessarily before $z$ does), the pair resumes moving 
along parallel trajectories inside the lower half-plane. Exactly the same line of calculation
as for Snell's law proves that the angles of incidence $\theta_1$, $\theta_2$ of these 
trajectories pre- and post-reflection satisfy, if $\mathrm{Im}\, \mu\neq0$,
\begin{equation}\label{eq:generalized reflection}
\sin\theta_1 \exp\Bigl(-\frac{{\rm Re}\, \mu}{{\rm Im}\, \mu}\, \theta_1 \Bigr) =
\sin\theta_2 \exp\Bigl(-\frac{{\rm Re}\, \mu}{{\rm Im}\, \mu}\, \theta_2 \Bigr).
\end{equation}
If $\mu$ is pure imaginary,  then (\ref{eq:generalized reflection}) reduces to
$\sin\theta_1 = \sin\theta_2$, from which we extract
\begin{equation}\label{eq:reflection}
\theta_1 = \pi - \theta_2 \qquad
{\rm if~~Im}\, \mu \neq 0\, , \> {\rm Re}\, \mu = 0\, . 
\end{equation}
This is the analogue of the law of reflection in optics, $\theta_1$ being the angle of
incidence, $\pi - \theta_2$ the angle of reflection.  Note that $s$ does not feature
in the results (\ref{eq:generalized reflection}) and (\ref{eq:reflection}).   The reason
is most readily appreciated in the case ${\rm Re}\, \mu = 0$: then reflection results
from `pivoting' of $z'(0) - z(0)$ to $z'(t) - z(t)$, and the relative magnitudes of $s_1, s_2$ 
determine where the pivot is (at the intersection of the dashed lines in Figure~\ref{fig:reflection}), 
but the net pivoting angle is the same regardless of where the pivot is.  

   If ${\rm Im}\, \mu = 0$, then no reflection occurs: $z'$ and $z$ repel each other along
a straight line, $z'$ running away in the upper half-plane and $z$ running away in the lower 
half-plane.  

   There arises a curious degeneracy when ${\rm Re}\, \mu = 0$ and the pair hits the mirror head-on, 
$\theta_1 = 0$.  The intuitive picture is that, barely into the upper half-plane,
both poles switch the sign of their strengths simultaneously, step
back barely into the lower half-plane, change the sign, etc.,
oscillating upper, lower, upper, lower, $\ldots$  The pair gets {\it trapped in the mirror}. This is reminiscent of work of Longuet-Higgins \cite{LH} on trapping of waves.

   We have detailed above two examples of the dynamics of a pole pair across a 
straight-line discontinuity of a piecewise constant `seabed' function.  In fact, these examples
are the basic ingredients from which we can construct almost all of the dynamics of {\it small 
pole pairs\/} on an arbitrary `seabed' $S$.  Here is the idea of the construction.

   Given an arbitrary $S$, approximate it by a piecewise constant function $S'$.  
The curves separating the level sets of $S'$ are generically smooth, so for pole pairs small compared
with the length-scale of the `seabed' geometry, each curve can be approximated 
locally by a straight line.  If $S'$ keeps the same sign on the two sides of such a line,
use the Snell calculation; if $S'$ changes the sign, use the reflection-off-the-mirror calculation.
In this approximation scheme, the only risk of misrepresenting the dynamics occurs 
if the pole pair hits a non-differentiable corner of one of the curves.  But non-differentiable
corners can be eliminated by perturbing the choice of $S$ or of $S'$, and anyway trajectories
that hit such corners are negligible in the space of all trajectories.  

   Now further optical analogues suggest themselves.  Just let loose an ensemble of
small pole pairs.  Since geometric optics is reducible to the laws (\ref{eq:Snell}) and 
(\ref{eq:reflection}), it is clear that our ensemble behaves like a pencil of light rays, although the 
behavior is more general and the dynamics richer if ${\rm Re}\, \mu \neq 0$.  The difference
of behavior is caused also by the pole-to-pole interaction between different pairs---light rays do
not interact among themselves---but this can be made small by making the pole pairs small.

\subsection{Leapfrogging and rainbow}

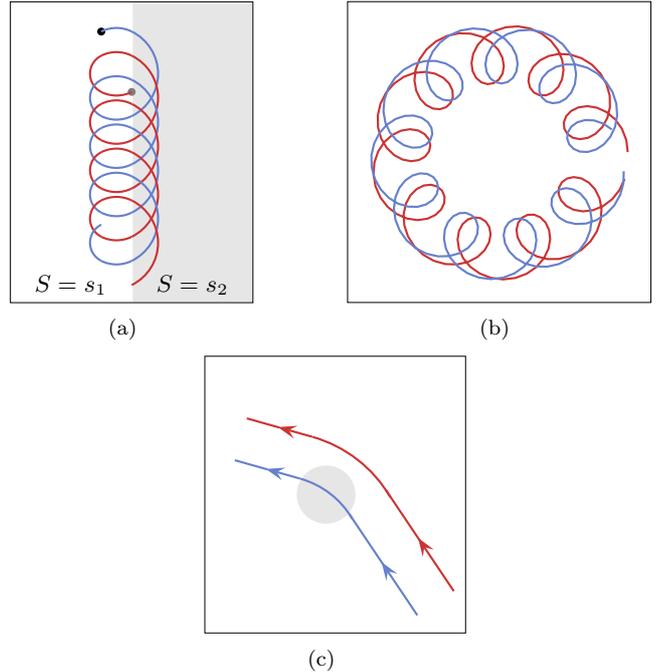
\begin{figure}[htb]
\begin{center}
\subfigure[]{
\psset{unit=0.8} %
 \begin{pspicture}(-2,-2)(2,3)
   \psframe[linestyle=none,fillstyle=solid,fillcolor=mygray](0,-2)(2,3)
   \psframe[linewidth=0.3pt](-2,-2)(2,3)
   \psdot[dotsize=3pt](-0.5,2.5) \psdot[dotsize=3pt,linecolor=gray](0,1.5)
  \multips{0}(0,0)(0,-0.8){4}{%
   \psarc[linecolor=myblue](-0.25,2){0.56}{63.43}{116.57}
   \psarc[linecolor=myred](-0.25,2){0.56}{-116.57}{-63.43}
   \psarc[linecolor=myblue](-0.35,1.8){0.784}{-63.43}{63.43}
   \psarc[linecolor=myred](-0.35,1.8){0.335}{116.57}{243.43}
   \psarc[linecolor=myred](-0.25,1.6){0.56}{63.43}{116.57}
   \psarc[linecolor=myblue](-0.25,1.6){0.56}{-116.57}{-63.43}
   \psarc[linecolor=myred](-0.35,1.4){0.784}{-63.43}{63.43}
   \psarc[linecolor=myblue](-0.35,1.4){0.335}{116.57}{243.43}
   }
   \rput(-1,-1.7){\small$S=s_1$}
   \rput(1,-1.7){\small$S=s_2$}
 \end{pspicture}}
\hspace{1cm}
\subfigure[]{%
\begin{pspicture}(-2,-2)(2,2)
  \psframe[linewidth=0.3pt](-2,-2)(2,2)\psset{unit=1.3}
  \parametricplot[plotpoints=200,linecolor=myred]{0}{350}{t cos t 10.23 mul cos 0.3 mul add t sin 0.3 t 10.23 mul sin mul add}
  \parametricplot[plotpoints=200,linecolor=myblue]{358}{705}{t cos t 10.23 mul cos 0.3 mul add t sin 0.3 t 10.23 mul sin mul add}
\end{pspicture}}
\hspace{1cm}
\subfigure[]{%
\psset{unit=0.8}
\begin{pspicture}(-2,-2.3)(2,2.3)
 \psframe[linewidth=0.3pt](-2,-2.3)(2.3,2.3)
 \pscircle[fillstyle=solid,fillcolor=mygray,linestyle=none](0,0){0.5}
{\psset{linecolor=myblue}
 \psline{->}(1.5,-2)(.91150, -1.1172) \psline(.94092, -1.1614)(.3819, -.3228)
 \psarc(-.8181, -1.1228){1.442}{30}{74}
 \psline{->}(-.42426, .2646)(-.99043, .42531) \psline(-.96212, .41727)(-1.5000, .56995)
}
{\psset{linecolor=myred}
 \psline{->}(2.1, -1.6)(1.5115, -.71725) \psline(1.5409, -.76138)(.9819, 0.0772)
 \psarc(-.8181, -1.1228){2.163}{30}{74}
 \psline{->}(-.22735, .9583)(-.79352, 1.1190) \psline(-.76521, 1.1110)(-1.3030, 1.2636)
}
\end{pspicture}}
  \caption{Leapfrogging and rainbow}\label{fig:leapfrogging}
\end{center}
\end{figure}

In vortex dynamics, {\it two} vortex pairs perform a motion called leapfrogging (e.g.\  \cite{MST} 
and references therein). With position-dependent strengths {\it one\/} pole pair can 
leapfrog all by itself.  
Let $S(z) = s_1$ in the left half-plane ${\rm Re}\, z < 0$ and $S(z) = s_2$ in the right
half-plane ${\rm Re}\, z \geqslant 0$, where $s_1, s_2 \in {\mathbb R}_+$.  A pair $z$, $z'$ with
strengths $\imagunit \, S(z)$, $\imagunit \, S(z')$ is initially at positions $\Delta z \in {\mathbb C}$ 
and $0$.  It is easy to check that this pair leapfrogs along piecewise circular paths, as in Figure~\ref{fig:leapfrogging}(a),
advancing by $\frac{s_2 - s_1}{s_2 + s_1}\imagunit \, |{\rm Im}\, \Delta z|$ per half-period.  
By distributing $S$ in a circular bump, like a plateau or a crater, we can also persuade a pair 
to leapfrog (quasi-) periodically as in Figure~\ref{fig:leapfrogging}(b).

   A pole pair of opposite $\mu = \pm \imagunit $ and of separation $d > 0$ approaches a disk of radius $r > 0$.  Let $S(z) = 1/r$ inside the disk and $S(z) =1$ outside.  Figure~\ref{fig:leapfrogging}(c) shows how the pair
gets refracted if one of the pair passes through the disk while the other misses it; the angle of refraction is
$2\, \arctan(1-r)/d$.  If both enter the disk, then the pair gets internally reflected a certain
number of times before it re-emerges, as in a rainbow.

In both Figures \ref{fig:leapfrogging}(b) and (c), the `seabed' function has circular symmetry.  The Noether-type theorem \ref{thm:Noether} applied to the rotational symmetry gives rise to a conservation law as described in (\ref{eq:rotational momentum}): namely conservation of
$$\psi(z,z') = \sigma(|z|^2) - \sigma(|z'|^2),$$
where $\sigma'(u)=S(u)$. In part (c), the conserved quantity is, when both vortices are outside the disk of radius $r$,
$$\psi(z,\,z') = |z|^2-|z'|^2.$$

\section{Numerical results for vortex pairs}

We show a sampling of striking numerical experiments for vortex pairs; that is, with $\mu'=-\mu$ both pure imaginary, and with variable `seabed' function.  Explicitly, we consider the system
\begin{equation}\label{eq:vortex pair}
\frac{d}{dt} z = \frac{\imagunit \, S(z')}{ \overline{z} - \overline{z'} }\, , \qquad
\frac{d}{dt} z' = - \frac{\imagunit \, S(z)}{ \overline{z'} - \overline{z} }\, .
\end{equation}
We consider $S(z)$ to be a combination of piecewise linear and constant functions.

The numerics were performed using \textsc{Maple}, with 
the default method and a step-size $\approx10^{-3}$, 
to integrate the system (\ref{eq:vortex pair}) with particular initial values and particular choices of the `seabed' function $S(z)$. The vortex pair is separated by about $5\times10^{-2}$ and the integration runs for the order of a few hundred time units.  There is much scope for investigations in a similar vein of other cases and with more poles. In the color figures, for each pair, the blue vortex has $\mu = +\imagunit $, the red $\mu' = -\imagunit $.  All
figures depict trajectories in the plane.  The pairs depicted in each figure are independent and not interacting with one another: they are merely different initial conditions.

Recall that in the usual vortex dynamics (with $S(z)= 1$ everywhere), an isolated pair moves along a straight line in a direction perpendicular to the segment that connects the vortices. We say the pair \emph{faces\/} this direction.
In each figure, the vortices in all pairs are separated by the same distance and, as reamrked before subsection \ref{sec:refraction}, this distance is a conserved quantity for arbitrary $S(z)$.

\subsection{Linear seabed}\label{sec:linearseabed}

\begin{figure}[htb]
\begin{center}
{\includegraphics[scale=0.4]{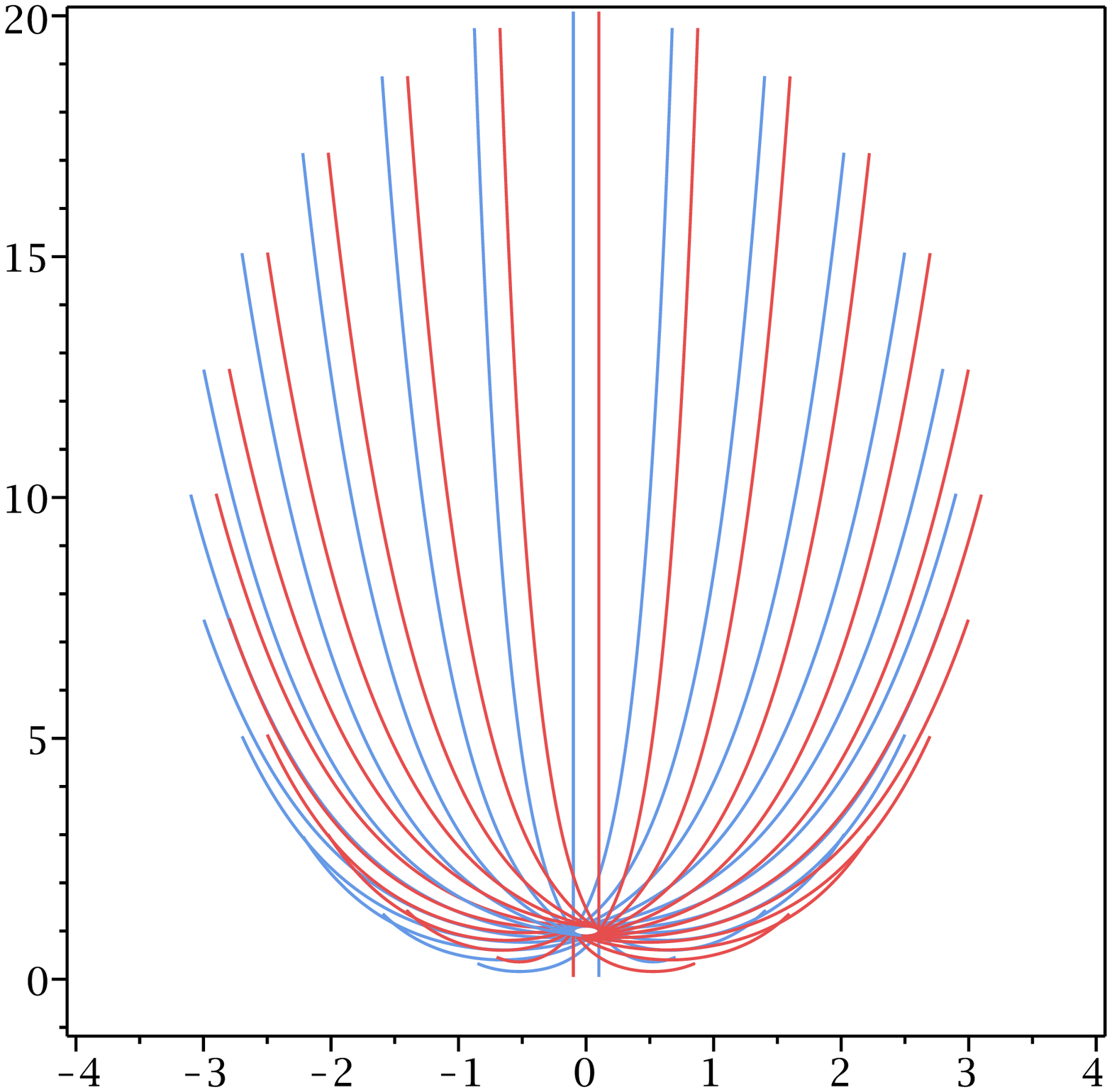}}
\caption{Linear Seabed}\label{fig:linear seabed}
\end{center}
\end{figure}

The `seabed' function is $S(z)=\Im z$. Initially the pairs are facing various  directions distributed uniformly around the circle, and they all start very near the point $z=\imagunit $, see Figure~\ref{fig:linear seabed}.  All the trajectories are bent upward, asymptotic to the vertical direction: a focusing effect. Notice the focusing  effect is stronger than depicted: the horizontal and vertical scales are different.  If the initial points are in the lower half plane, where $S(z)<0$, they move downward in the same way. 

(In general, replacing $S(z)$ by $-S(z)$ is equivalent to reversing vorticity---so interchanging the pair, or to reversing time.)


In this example, and those of the following subsection, the `horizontal' (real) translation is a symmetry of $S$ so the Noether-type theorem  \ref{thm:Noether} gives rise to a conservation law.  The conserved `momentum' is given by (\ref{eq:translational momentum}), with $\sigma'(y)=y$. Thus one can take $\sigma(y)=\frac12 y^2$, so that $\psi(z,z') = \frac12(y^2-y'^2)$
is a conserved quantity.

\subsection{Piecewise linear seabeds}
\label{sec:piecewise linear seabed} 

In this series of experiments (a)--(d), $S(z)$ is a piecewise linear, continuous function of $\Im z$. 

\bigskip

\noindent(a) Here we take 
$$
S(z) = \begin{cases}
1 & \text{if~~ } \Im z < 0\, ,\\
1+\Im z & \text{if~~ } 0 \leqslant \Im z < 1\, ,\\
2 & \text{if~~ }\Im z  \geqslant 1
\end{cases}
$$
whose graph as a function of $\Im z$ is
$\psset{unit=0.4}\begin{pspicture}(-2,1)(3,2)
\psline(-2,1)(0,1)(1,2)(3,2)
\end{pspicture}$~. 

\begin{figure}[htb]
\begin{center}
\includegraphics[scale=0.35]{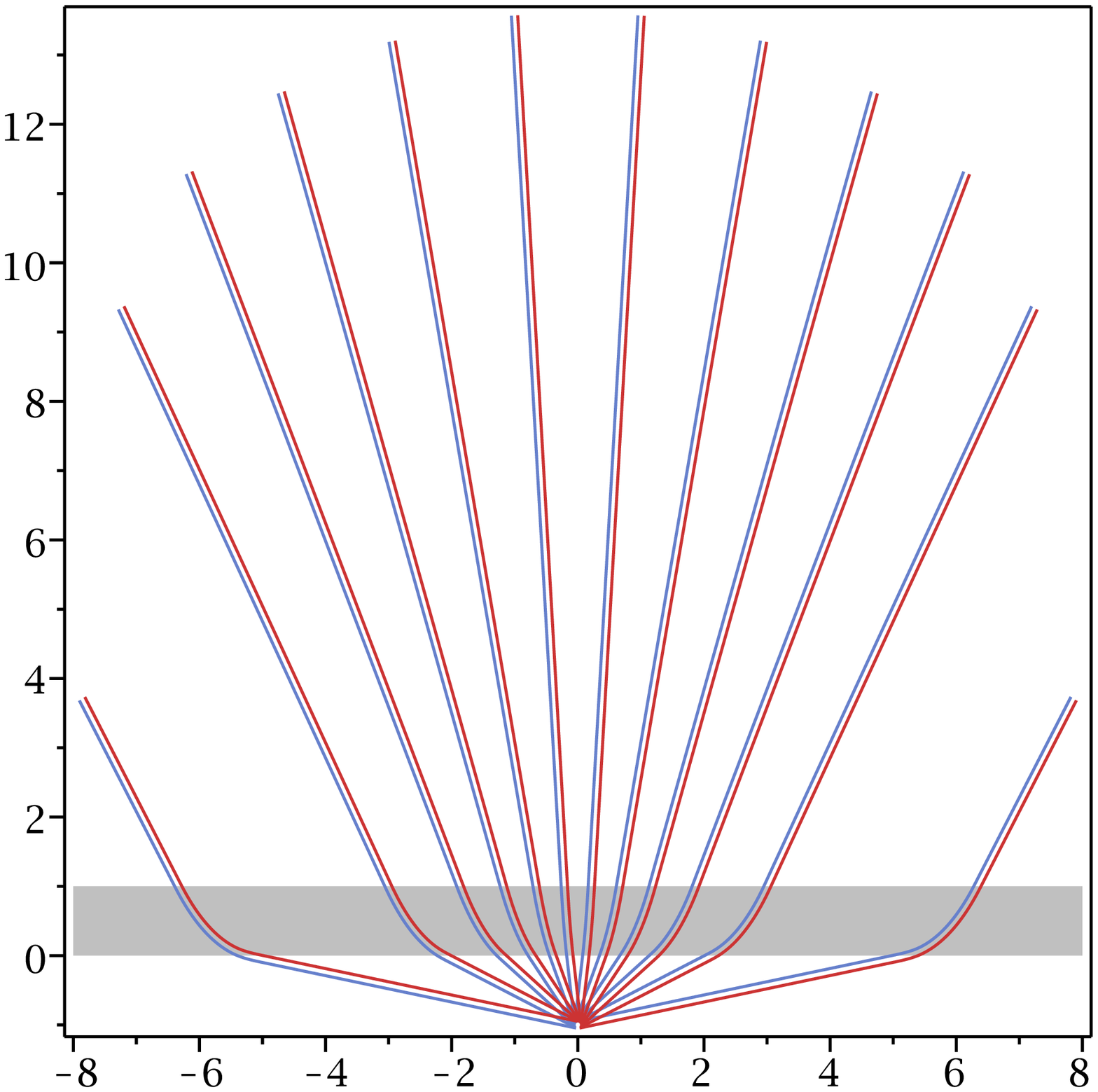}
\caption{Sloped step}\label{fig:spread-refraction}
\end{center}
\end{figure}

Initially the pairs face various angular directions and they start near the point $z=-\imagunit $.  Pairs that initially face downward will continue to move downward as $S$ there is flat, and these are not depicted.  On the other hand, if they have an initial upward component, then they enter the sloping part of $S$ (shown in grey in Figure~\ref{fig:spread-refraction}), incur a certain measure of focusing as described in subsection \ref{sec:linearseabed}, and then emerge on the flat part, moving thereafter along straight lines.  We see an overall refraction pattern much as described in subsection \ref{sec:refraction} above.

\begin{figure}[htb]
\begin{center}
\includegraphics[scale=0.35]{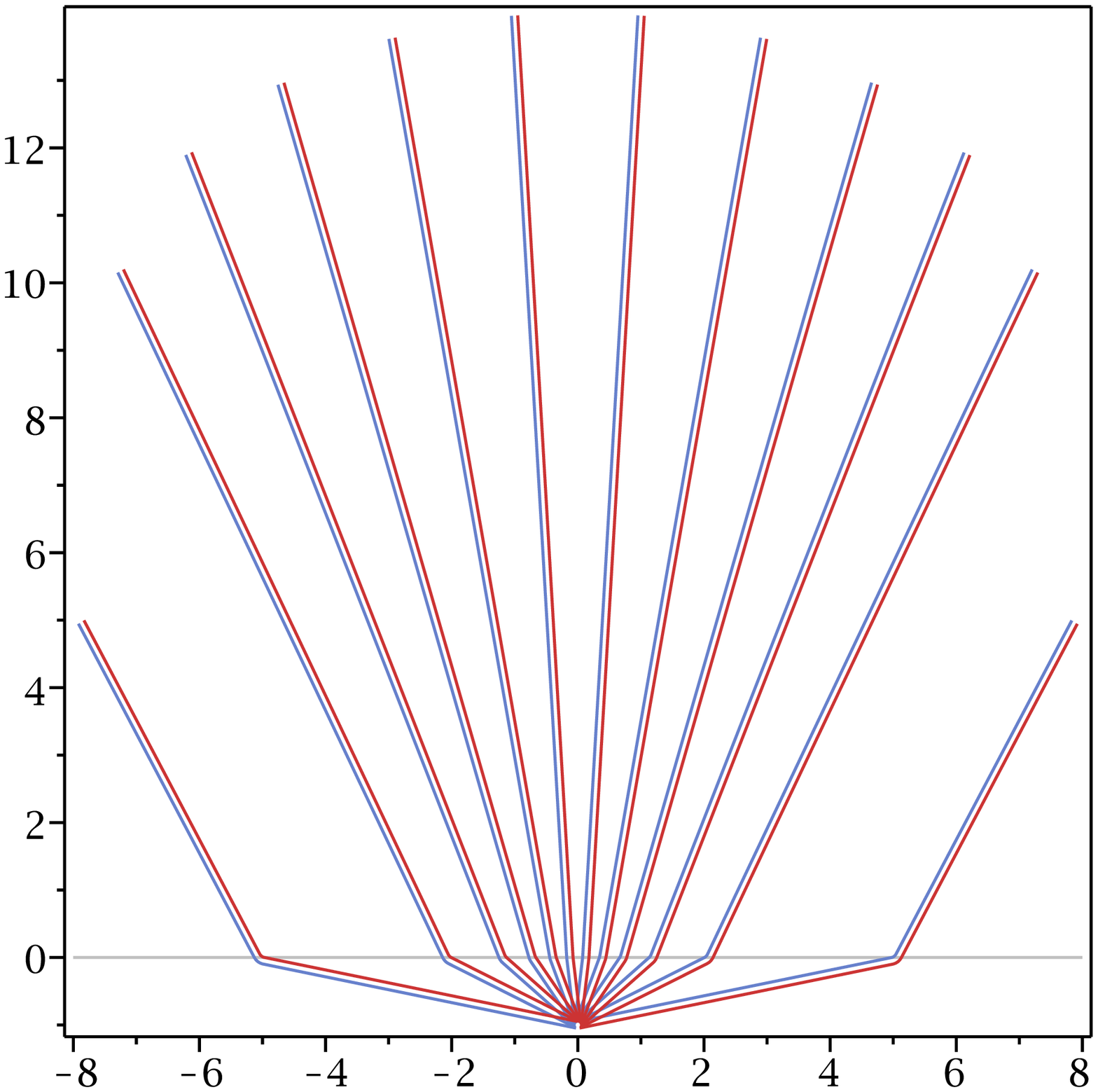}
\caption{Abrupt step}\label{fig:abrupt step}
\end{center}
\end{figure}

The results can be compared with the analogous figure where the change of `level' is abrupt, which is precisely the setting of subsection \ref{sec:refraction}; see Figure~\ref{fig:abrupt step}, where the `seabed' changes from 1 to 2 as $\Im(z)$ passes through 0. 

\bigskip

\begin{figure}[htb]
\begin{center}
\includegraphics[scale=0.3]{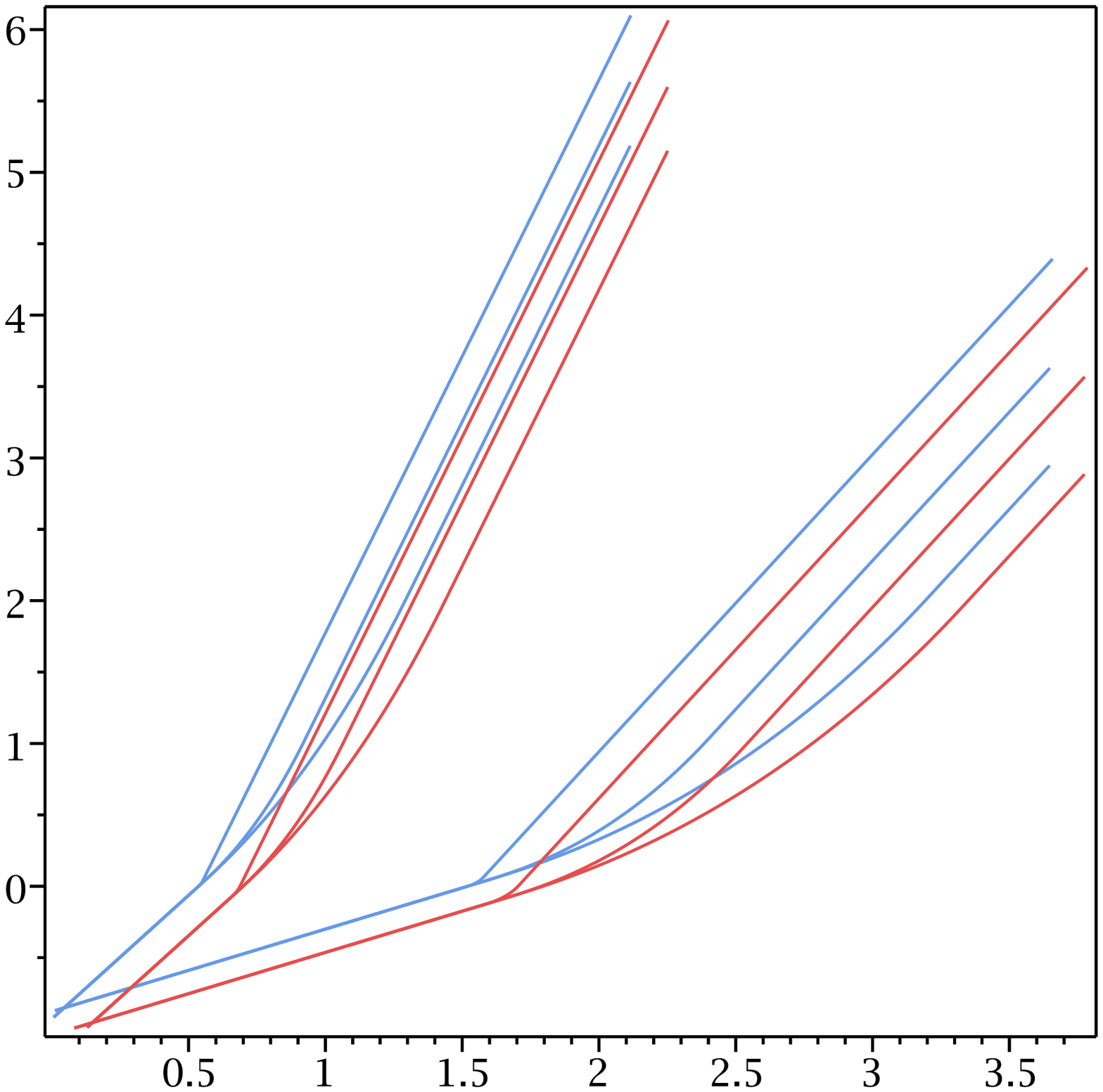}
\caption{Comparison of interpolating slopes} \label{fig:different slopes}
\end{center}
\end{figure}

\noindent(b) We next examine how the sloping part of the type of `seabed' used in (a) affects the angle of refraction.
In  Figure~\ref{fig:different slopes}, two vortex pairs initially face the angular directions $\pi/3$ and $\pi/6$ as measured from the 
real axis, and start near $z = -1$.  Each pair is then subjected to three different `seabed' functions.

The trajectory of each pair therefore splits into three offshoots: the upper offshoot is subjected to a step function $S(z)$ with a discontinuity across $\Im z = 0$ (as above), whose graph as a function of $\Im z$ is
$$
\psset{unit=0.4}\begin{pspicture}(-2,1)(3,2)
\psline(-2,1)(0,1)(0,2)(3,2)
\end{pspicture}
$$ 
or in formulae
$$
S(z) = \begin{cases}
1 & \text{if~~ } \Im z < 0\, ,\\
2 & \text{if~~ } \Im z \geqslant 0\, .
\end{cases}
$$
The middle offshoot is subjected to the same `seabed' as in the previous experiment (a), whose graph was
$$
\psset{unit=0.4}\begin{pspicture}(-2,1)(3,2)
\psline(-2,1)(0,1)(1,2)(3,2)
\end{pspicture}
$$
The lower offshoot is subjected to a similar `seabed' but whose sloping part has gradient $1/2$, stretched over the wider strip $0\leqslant \Im z < 2$ :
$$
\psset{unit=0.4}\begin{pspicture}(-2,1)(3,2)
\psline(-2,1)(0,1)(2,2)(3,2)
\end{pspicture}
$$
In all three the function increases from 1 to 2 in value.

The fact that for each pair the three offshoots end up parallel means that the net angle of refraction is independent of the type of `seabed' to which an offshoot was subjected.  Similar experiments suggest the same is true if the `seabed' is nonlinear on the interpolating strip. In the figure, it is the step function that deviates abruptly along the trajectories: the greater the slope, the faster the deviation. 


\bigskip

\noindent(c) Now consider the `seabed' function
$$
S(z) = \begin{cases}
-1 & \text{if~~ } \Im z < 0\, ,\\
-1+\Im z & \text{if~~ } 0 \leqslant \Im z < 2\, ,\\
1 & \text{if~~ }\Im z  \geqslant 2
\end{cases}
$$
whose graph as a function of $\Im z$ is
$\raisebox{-3pt}{\psset{unit=0.4}\begin{pspicture}(-2,1)(3,2)
\psline(-2,1)(0,1)(1,2)(3,2)
\psline[linecolor=gray,arrowsize=0.4]{->}(-2,1.5)(3,1.5)
\end{pspicture}}$~, 
which is like experiment (a) except that this time $S$ is negative on the lower half-plane.

\begin{figure}[htb]
\begin{center}
\includegraphics[scale=0.45]{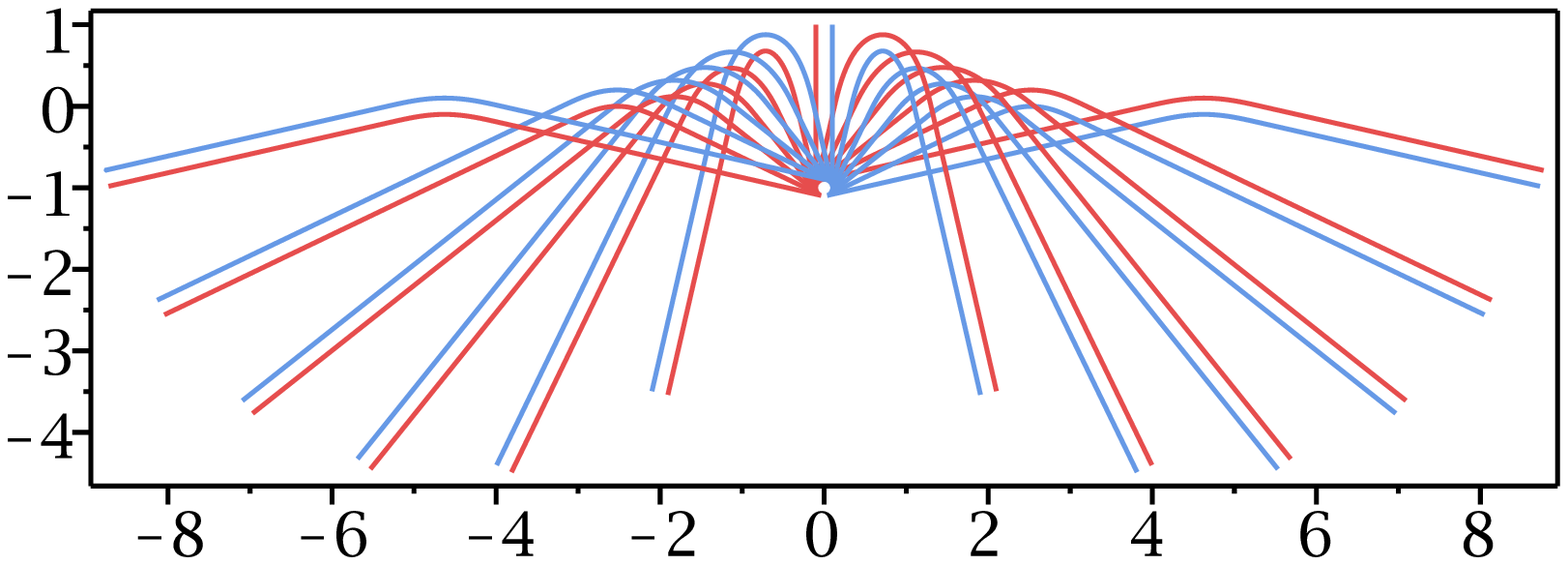}
\end{center}
\vspace{-5mm}
\caption{Experiment (c)}\label{fig:spider}
\end{figure}

This may be viewed as a continuous version of the mirror of subsection \ref{sec:reflection}; here too, the pairs undergo overall reflection, but closer inspection reveals that the change of direction is less sudden.  The penetration of the pair into the strip where $S$ slopes depends on its angle of departure. See Figure \ref{fig:spider}.

\bigskip

\begin{figure}[htb]
\begin{center}
\includegraphics[scale=0.35]{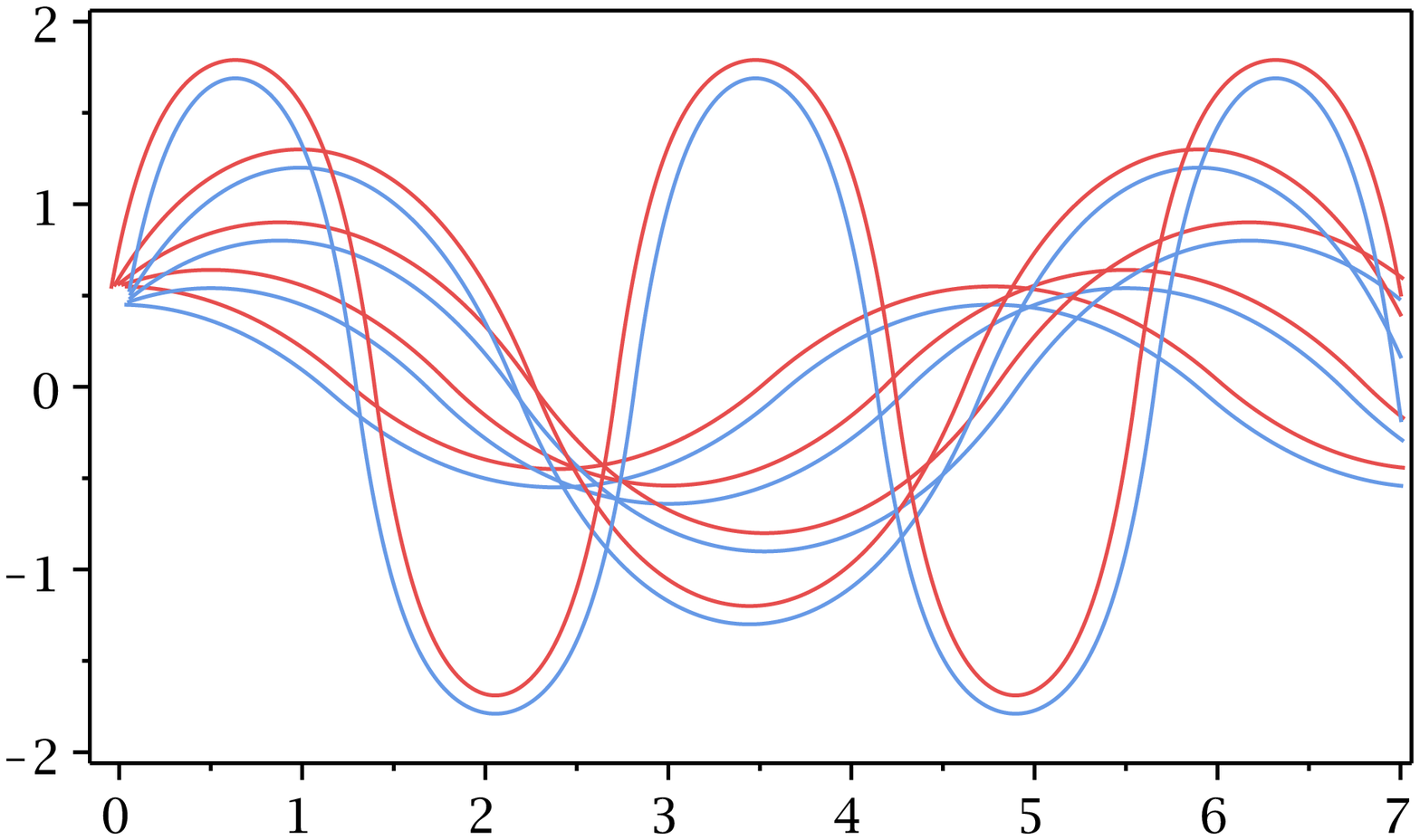}
\end{center}
\vspace{-5mm}
\caption{Vortex pairs trapped in mirrors} \label{fig:trapped}
\end{figure}

\noindent(d) 
Placing together two mirrors from (c) face to face, we take 
$$
S(z)=|\Im z| - 2
$$ 
whose graph as a function of $\Im z$ is
\raisebox{-10pt}{$\psset{unit=0.2}\begin{pspicture}(-4,-2)(4.2,2)
\psline(-4,2)(0,-2)(4,2)
\psline[linecolor=gray,arrowsize=0.7]{->}(-4,0)(4,0)
\end{pspicture}$}~. 
On the plane the graph of $S$ looks like a trough.  
Pairs released from near $z = 0$ where $S(z) < 0$ move in a kind of smoothed zigzag, drifting along the real axis.  Their trajectories are confined within a certain strip near the bottom of the trough, as depicted in Figure~\ref{fig:trapped}.  The width of the confinement strip depends on the direction the pair initially faces.  It depends little on the strength or the size of the pair, for the same reason that these two parameters do not feature in the law of reflection (\ref{eq:generalized reflection}) and (\ref{eq:reflection}).

\begin{figure}[ht]
\begin{center}
\includegraphics[scale=0.3]{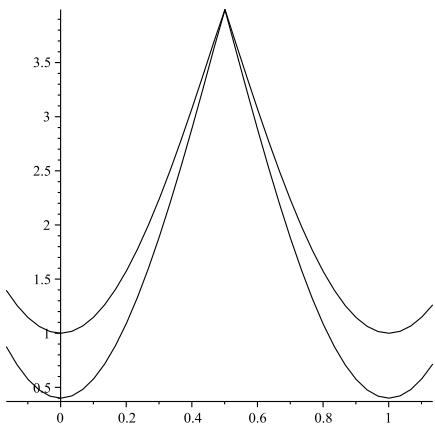}
\end{center}
\vspace{-5mm}
\caption{Amplitudes of trapped vortex pairs}\label{fig:amplitude}
\end{figure}

The amplitude of the smooth zigzags is shown in Figure~\ref{fig:amplitude}, as a 
function of the angle of departure $\theta$ with the horizontal (the abscissa is $\theta/\pi$). The two curves are for different initial values of $\Im z$; the lower one for $\Im z=0.2$, the upper for $\Im z=0.5$ (where $z$ is the mid-point of the initial position of the pair). As can be seen in the previous figure, zigzags take place in the strip $-2 < \Im z < 2$ where $S(z) < 0$. 
The maximum of the curves is attained by pairs facing almost vertically upward.  (If the pair faces exactly upward, then it tends to $\Im z = 2$ where $S(z) =0$ in infinite time.) The minimum occurs where the pair faces horizontally, and the motion has the initial point as an extremum.

Confinement likewise occurs with a step function version of this, whose graph as a function of $\Im z$ is
\raisebox{-12pt}{$\psset{unit=0.4}\begin{pspicture}(-2,-1)(2.2,1)
\psline(-2,1)(-1,1)(-1,-1)(1,-1)(1,1)(2,1)
\psline[linecolor=gray]{->}(-2,0)(2,0)
\end{pspicture}$}~. 
The trajectories, which resemble those in the trough above, are readily derived by iterating formulae for the mirror in subsection \ref{sec:reflection}.

\subsection{Caustic from a circular mirror}
\label{sec:caustic}

\begin{figure}[thb]
\begin{center}
\includegraphics[scale=0.4]{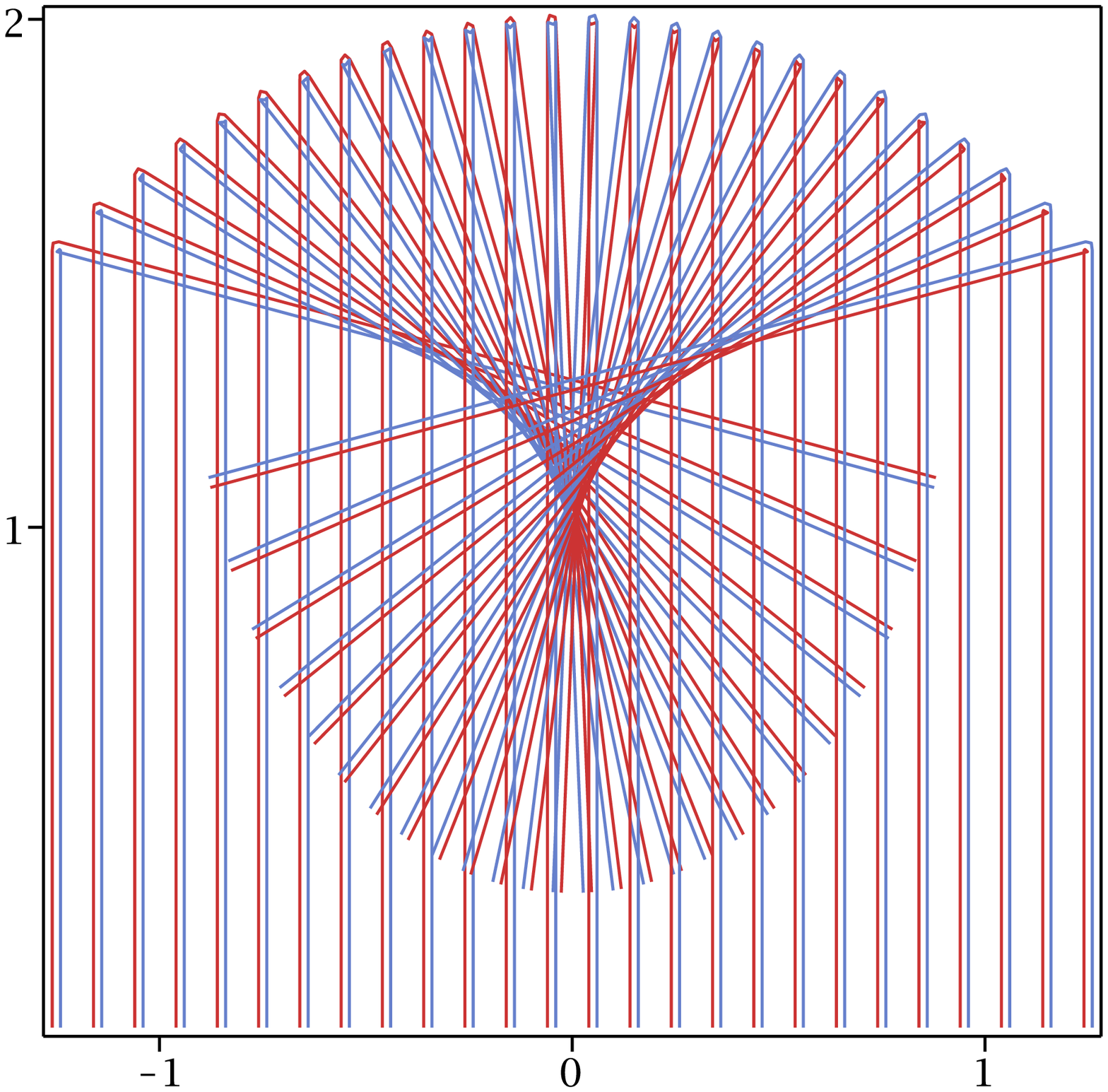}
\end{center}
\vspace{-5mm}
\caption{Caustic from a circular mirror}\label{fig:caustic}
\end{figure}

We take a mirror in the form of a circular arc (of radius 2 in the figure).
The `seabed' function is equal to $-1$ below the arc, $+1$ above the arc.  An ensemble of pairs is released, each pair initially
facing vertically upward along the bottom of the figure (below the mirror).  They strike the mirror at the top, and are reflected much 
as described in subsection \ref{sec:reflection}.  Together they form a caustic, familiar from when we observe the pattern of light on the surface of tea inside a tea cup.


\paragraph{Acknowledgements} 
 We would like to thank the referees for comments improving the original manuscript and pointing out a number of relevant references.


\end{document}